# Markovianity in space and time[*]

M. N. M. van Lieshout[1]

*Centre for Mathematics and Computer Science, Amsterdam*

**Abstract.** Markov chains in time, such as simple random walks, are at the heart of probability. In space, due to the absence of an obvious definition of past and future, a range of definitions of Markovianity have been proposed. In this paper, after a brief review, we introduce a new concept of Markovianity that aims to combine spatial and temporal conditional independence.

## 1. From Markov chain to Markov point process, and beyond

This paper is devoted to the fundamental concept of *Markovianity*. Although its precise definition depends on the context, common ingredients are conditional independence and factorisation formulae that allow to break up complex, or high dimensional, probabilities into manageable, lower dimensional components. Thus, computations can be greatly simplified, sometimes to the point that a detailed probabilistic analysis is possible. If that cannot be done, feasible, efficient simulation algorithms that exploit the relatively simple building blocks may usually be designed instead.

### 1.1. Markov chains

The family of *Markov chains* is one of the most fundamental and intensively studied classes of stochastic processes, see e.g. [2]. If we restrict ourselves to a bounded time horizon, say $0, 1, \ldots, N$ for some $N \in \mathbb{N}$, a discrete Markov chain is a random vector $(X_0, \ldots, X_N)$ with values in some denumerable set $L \neq \emptyset$ for which

$$\mathbb{P}(X_i = x_i \mid X_0 = x_0; \ldots; X_{i-1} = x_{i-1}) = \mathbb{P}(X_i = x_i \mid X_{i-1} = x_{i-1}) \qquad (1)$$

for all $1 \leq i \leq N$ and all $x_j \in L$, $j = 0, \ldots, i$. In words, the Markov property (1) means that the probabilistic behaviour of the chain at some time $i$ given knowledge of its complete past depends only on its state at the immediate past $i-1$, regardless of how it got to $x_{i-1}$. The right hand side of equation (1) is referred to as the *transition probability* at time $i$ from $x_{i-1}$ to $x_i$. If $\mathbb{P}(X_i = x_i \mid X_{i-1} = x_{i-1}) = p(x_{i-1}, x_i)$ does not depend on $i$, the Markov chain is said to be *stationary*.

By the product rule, the joint distribution can be factorised as

$$\mathbb{P}(X_0 = x_0; \ldots; X_N = x_N) = \mathbb{P}(X_0 = x_0) \prod_{i=1}^{N} p(x_{i-1}, x_i). \qquad (2)$$

---

[*]This research is supported by the Technology Foundation STW, applied science division of NWO, and the technology programme of the Ministry of Economic Affairs (project CWI.6156 'Markov sequential point processes for image analysis and statistical physics').

[1]CWI, P.O. Box 94079, 1090 GB Amsterdam, The Netherlands, e-mail: colette@cwi.nl

Consequently,

$$\mathbb{P}(X_i = x_i \mid X_j = x_j, j \neq i) = \mathbb{P}(X_i = x_i \mid X_{i-1} = x_{i-1}, X_{i+1} = x_{i+1}) \quad (3)$$

for all $i = 1, \ldots, N-1$, and all $x_j \in L$, $j = 0, \ldots, N$ (with obvious modifications at the extremes $i = 0, N$ of the time interval of interest). Thus, the conditional distribution of the state at a single point in time depends only on the states at the immediate past and future. To tie in with a more general concept of Markovianity to be discussed below, the time slots $i-1$, $i+1$ (if within the finite horizon) may be called *neighbours* of $i$, $i \in \{0, \ldots, N\}$.

As a simple example, consider the celebrated *simple random walk* on the two dimensional lattice $L = \mathbb{Z}^2$. The dynamics are as follows: A particle currently sitting at site $x \in L$ moves to each of the neighbouring sites with equal probability. More precisely, $p(x, y) = 1/4$ if $x$ and $y$ are horizontally or vertically adjacent in $L$, that is, $||x - y|| = 1$, and zero otherwise [25].

### 1.2. Markov random fields

The concept of Markovianity plays an important role in space as well as in time, especially in image analysis and statistical physics. Since, in contrast to the setting in the previous subsection, space does not allow a natural order, a symmetric, reflexive *neighbourhood* relation $\sim$ must be defined on the domain of definition, which we assume to be a finite set $I$. The generic example is a rectangular grid $I \subseteq \mathbb{Z}^2$ with $i \sim j$ if and only if $||i - j|| \leq 1$. Then, a discrete *Markov random field* with respect to $\sim$ is a random vector $X = (X_i)_{i \in I}$, with values of the components $X_i$ in some denumerable set $L \neq \emptyset$, that satisfies the *local Markov property*

$$\mathbb{P}(X_i = x_i \mid X_j = x_j, j \neq i) = \mathbb{P}(X_i = x_i \mid X_j = x_j, j \in \partial(i))$$

where $\partial(i) = \{j \in I \setminus \{i\} : j \sim i\}$ is the set of neighbours of $i$. The expression should be compared to (3); indeed a Markov chain is a Markov random field on $I = \{0, 1, \ldots, N\}$ with $i \sim j$ if and only if $|i - j| \leq 1$, $i, j \in I$.

To state the analogue of (2), define a *clique* to be a set of sites $C \subseteq I$ such that for all $i, j \in C$ these sites are neighbours, that is $i \sim j$. By default, the empty set is a clique. If we assume for simplicity that the joint probability mass function of the random vector $X$ has no zeroes, it defines a Markov random field if and only if it can be factorised as

$$\mathbb{P}(X_0 = x_0; \ldots; X_N = x_N) = \prod_{\text{cliques } C} \varphi_C(x_j, j \in C) \quad (4)$$

for some *clique interaction functions* $\varphi_C(\cdot) > 0$. Hence, the clique interaction functions are the spatial analogues of the transition probabilities of a Markov chain.

A two-dimensional example is the *Ising model* for spins in a magnetic field [15]. In this model, each node site of a finite subset of the lattice is assigned a spin value from the set $L = \{-1, 1\}$. Interaction occurs between spins at horizontally or vertically adjacent sites, so that cliques consist of at most two points. A particular spin value may be preferred due to the presence of an external magnetic field. More precisely, the Ising model is defined by (4) with $\varphi_{\{i\}}(x_i) = e^{\alpha x_i}$ for singletons, $\varphi_{\{i,j\}}(x_i, x_j) = e^{\beta x_i x_j}$, and $\varphi_\emptyset$ set to the unique constant for which the right hand side of (4) is a probability mass function. The constant $\alpha$ reflects the external magnetic field and influences the frequencies of the two spin types. For $\beta > 0$, neighbouring sites tend to agree in spin, for $\beta < 0$, they tend to have different spins. For further details, the reader is referred to [11].

## *1.3. Markov point processes*

The next step is to leave discrete domains and move into Euclidean space. Thus, let $D \subset \mathbb{R}^d$ be a disc, rectangle or other bounded set of positive volume. A realisation of a *Markov point process* on $D$ is a subset $\mathbf{x} = \{x_1, \ldots, x_n\}$ of $D$ with random, finite $n \geq 0$. Again define a reflexive, symmetric neighbourhood relation $\sim$, for instance by $x \sim y$ if and only if $||x - y|| \leq R$ for some $R > 0$. Note that each point has uncountably many potential neighbours, in contrast to the discrete case. Moreover, the probability of finding a point at any particular location $a \in D$ will usually be 0, so densities rather than probability mass functions are needed. A suitable dominating measure is the homogeneous Poisson process (see e.g. [17]), because of its lack of spatial interaction. Indeed, given such a Poisson process places $n$ points in $D$, the locations are i.i.d. and uniformly distributed over $D$. If desirable, one could easily attach *marks* to the points, for example a measurement or type label of an object located at the point. Doing so, the domain becomes $D \times M$, where $M$ is the mark set, and each $x_i$ can be written as $(d_i, m_i)$.

In this set-up, a *Markov point process* $X$ is defined by a density $f(\cdot)$ that is hereditary in the sense that $f(\mathbf{x}) > 0 \Rightarrow f(\mathbf{y}) > 0$ for all configurations $\mathbf{y} \subseteq \mathbf{x}$ and satisfies [18, 26] the *local Markov condition*, that is, whenever $f(\mathbf{x}) > 0$, $u \notin \mathbf{x}$, the ratio

$$\lambda(u \mid \mathbf{x}) := \frac{f(\mathbf{x} \cup \{u\})}{f(\mathbf{x})} \tag{5}$$

depends only on $u$ and $\{x_i : u \sim x_i\}$. The function $\lambda(\cdot \mid \cdot)$ is usually referred to as *conditional intensity*. For the null set where $f(\mathbf{x}) = 0$, $\lambda(\cdot \mid \mathbf{x})$ may be defined arbitrarily. Note that since realisations $\mathbf{x}$ of $X$ are *sets* of (marked) points, $f(\cdot)$ is symmetric, i.e. invariant under permutations of the $x_i$ that constitute $\mathbf{x}$. If

$$f(\mathbf{x} \cup \{u\}) \leq \beta f(\mathbf{x})$$

for some $\beta > 0$ uniformly in $\mathbf{x}$ and $u$, the density $f(\cdot)$ is said to be *locally stable*. Note that local stability implies that $f(\cdot)$ is hereditary, and that the conditional intensity (5) is uniformly bounded in both its arguments.

A factorisation is provided by the *Hammersley–Clifford theorem* which states [1, 26] that a marked point process with density $f(\cdot)$ is Markov if and only if

$$f(\mathbf{x}) = \prod_{\text{cliques } \mathbf{y} \subseteq \mathbf{x}} \varphi(\mathbf{y}) \tag{6}$$

for some non-negative, measurable interaction functions $\varphi(\cdot)$. The resemblance to (4) is obvious.

An example of a locally stable Markov point process is the *hard core model*, a homogeneous Poisson process conditioned to contain no pair of points that are closer than $R$ to one another. For the hard core model, $\varphi(\{x, y\}) = \mathbf{1}\{||x - y|| > R\}$, $\varphi(\mathbf{y}) \equiv 1$ for singletons and configurations $\mathbf{y}$ containing three or more points. The interaction function for the empty set acts as normalising constant to make sure that $f(\cdot)$ integrates to unity.

Although Markov point processes are useful modelling tools [18], the assumption of permutation invariance may be too restrictive. For instance in image interpretation, occlusion of one object by another may be dealt with by an ordering of the objects in terms of proximity to the camera [20], or in scene modelling various kinds of alignment of the parts that make up an objects may be modelled by means of non-symmetric neighbour relations [22]. Non-symmetric densities also arise naturally when local scale and intensity is taken into account by transforming the conditional intensity of a homogeneous template process, see [14].

In stark contrast to the vast literature on spatial point processes, sequential patterns have not been studied much. Indeed, an analogue to the theory of Markov point processes does not exist in the non-symmetrical setting. The present paper fills this gap.

## 2. Definitions and notation

This paper is concerned with *finite sequential spatial processes*. Realisations of such processes at least consist of a finite sequence

$$\vec{\mathbf{x}} = (x_1, \ldots, x_n), \qquad n \in \mathbb{N}_0$$

of points in some bounded subset $D$ of the plane. Additionally, to each point $x_i$, a mark $m_i$ in some complete, separable metric space, say $M$, may be attached. The mark may be a discrete type label, a real valued measurement taken at each location, or a vector of shape parameters to represent an object located at $x_i$. Thus, we may write

$$\vec{y} = (y_1, \ldots, y_n) = ((x_1, m_1), \ldots, (x_n, m_n)), \qquad n \in \mathbb{N}_0,$$

for the configuration of marked points, and shall denote the family of all such configurations by $N^{\mathrm{f}}$.

As an aside, the plane $\mathbb{R}^2$ may be replaced by $\mathbb{R}^d$ or any other complete separable metric space, equipped with the Borel $\sigma$-algebra and a finite diffuse Borel measure.

The distribution of a finite sequential spatial process may be defined as follows. Given a finite diffuse Borel measure $\mu(\cdot)$ on $(D, \mathcal{B}_D)$ so that $\mu(D) > 0$, usually Lebesgue measure, and a mark probability measure $\mu_M(\cdot)$ on $M$ equipped with its Borel $\sigma$-algebra $\mathcal{M}$, specify

(1) a probability mass function $q_n$, $n \in \mathbb{N}_0$, for the number of points in $D$;
(2) for each $n$, a Borel measurable and $(\mu \times \mu_M)^n$-integrable joint probability density $p_n(y_1, \ldots, y_n)$ for the sequence of marked points $y_1, \ldots, y_n \in D \times M$, given it has length $n$.

Alternatively, a probability density $f(\cdot)$ may be specified directly on $N^{\mathrm{f}} = \cup_{n=0}^{\infty}(D \times M)^n$, the space of finite point configurations in $D$ with marks in $M$, with respect to the reference measure $\nu(\cdot)$ defined by $\nu(F)$ equal to

$$\sum_{n=0}^{\infty} \frac{e^{-\mu(D)}}{n!} \int_{D \times M} \cdots \int_{D \times M} \mathbf{1}\left\{(y_1, \ldots, y_n) \in F\right\} d\mu \times \mu_M(y_1) \cdots d\mu \times \mu_M(y_n)$$

for $F$ in the $\sigma$-algebra on finite marked point sequences generated by the Borel product $\sigma$-fields on $(D \times M)^n$. In words, $\nu(\cdot)$ corresponds to a random sequence of Poisson length with independent components distributed according to the normalised reference measure $\mu(\cdot) \mu_M(\cdot)/\mu(D)$.

It is readily observed that $q_0 = \exp(-\mu(D)) f(\emptyset)$, and

$$q_n = \frac{e^{-\mu(D)}}{n!} \int_{D \times M} \cdots \int_{D \times M} f(y_1, \ldots, y_n) \, d\mu \times \mu_M(y_1) \cdots d\mu \times \mu_M(y_n);$$

$$p_n(\vec{y}) = \frac{e^{-\mu(D)}}{n! \, q_n} f(\vec{y})$$

for each $n \in \mathbb{N}$ and $\vec{\mathbf{y}} \in (D \times M)^n$.

Reversely, if the length $n(\vec{\mathbf{y}})$ of $\vec{\mathbf{y}}$ is $n$,

$$f(\vec{\mathbf{y}}) = e^{\mu(D)} \, n! \, q_n \, p_n(\vec{\mathbf{y}}).$$

Note that neither $f(\cdot)$ nor the $p_n(\cdot, \ldots, \cdot)$ are required to be symmetric [3, Ch. 5].

We are now ready to define and analyse a Markov concept for random sequences in the plane. To do so, we begin by defining a *sequential conditional intensity*

$$\lambda_i((x, m) \mid \vec{\mathbf{y}}) := \frac{f(s_i(\vec{\mathbf{y}}, (x, m)))}{(n+1) f(\vec{\mathbf{y}})} \quad \text{if } f(\vec{\mathbf{y}}) > 0 \tag{7}$$

for inserting $(x, m) \notin \vec{\mathbf{y}}$ at position $i \in \{1, \ldots, n+1\}$ of $\vec{\mathbf{y}} = (y_1, \ldots, y_n)$. Here $s_i(\vec{\mathbf{y}}, (x, m)) = (y_1, \ldots, y_{i-1}, (x, m), y_i, \ldots, y_n)$. On the null set $\{\vec{\mathbf{y}} \in N^{\mathrm{f}} : f(\vec{\mathbf{y}}) = 0\}$, the sequential conditional intensity may be defined arbitrarily. The overall conditional probability of finding a marked point at $du = d\mu \times \mu_M(u)$ in any position in the vector given that the remainder of the sequence equals $\vec{\mathbf{y}}$ is given by

$$\sum_{i=1}^{n+1} \lambda_i(du \mid \vec{\mathbf{y}}).$$

The expression should be compared to its classic counterpart (5). As for Markov chains, we are mostly interested in $\lambda_{n+1}(\cdot \mid \vec{\mathbf{y}})$, but all $\lambda_i(\cdot \mid \cdot)$ are needed for the reversibility of the dynamic representation to be considered in Section 4.

Note that provided $f(\cdot)$ is hereditary in the sense that $f(\vec{\mathbf{y}}) > 0$ implies $f(\vec{\mathbf{z}}) > 0$ for all subsequences $\vec{\mathbf{z}}$ of $\vec{\mathbf{y}}$, then

$$f(y_1, \ldots, y_n) \propto n! \prod_{i=1}^{n} \lambda_i(y_i \mid \vec{\mathbf{y}}_{<i}) \tag{8}$$

where $\vec{\mathbf{y}}_{<i} = (y_1, \ldots, y_{i-1})$. This observation implies an alternative way of defining a sequential spatial process: by its conditional intensity. Some care has to be taken to make sure that the resulting density (8) is integrable with respect to $\nu(\cdot)$. A sufficient condition is that $\lambda_i(y_i \mid \vec{\mathbf{y}}_{<i}) \leq \beta/i$ for all choices of $i$, $y_i$, and $\vec{\mathbf{y}}_{<i}$. It is important to note, though, that such an assumption does not imply that $f(s_j(\vec{\mathbf{y}}, u))/f(\vec{\mathbf{y}})$ is uniformly bounded in $u$ and $\vec{\mathbf{y}}$ for $j \neq n(\vec{\mathbf{y}}) + 1$. The latter stronger assumption, which may be referred to as *local stability* in analogy with classic marked point processes, will be required in Section 4.

In order to generalise the concept of Ripley–Kelly Markovianity [18, 26] to random sequences, suppose a reflexive relation $\sim$ on $D \times M$ is given with the purpose of formalising the local interactions. In contrast to the point process context, we do not require $\sim$ to be symmetric. If $y \sim z$, the marked point $z$ is said to be a *directed neighbour* of $y$.

**Definition 1.** A sequential spatial process $Y$ on a bounded Borel set $D \subseteq \mathbb{R}^2$ with marks in $M$ defined by its density $f(\cdot)$ with respect to $\nu(\cdot)$ is said to be Markov with respect to the relation $\sim$ on $D \times M$ if

- $f(\vec{\mathbf{y}}) > 0$ implies $f(\vec{\mathbf{z}}) > 0$ for all subsequences $\vec{\mathbf{z}}$ of $\vec{\mathbf{y}}$;
- for all sequences $\vec{\mathbf{y}}$ for which $f(\vec{\mathbf{y}}) > 0$, the ratio $f((\vec{\mathbf{y}}, u))/f(\vec{\mathbf{y}})$ depends only on $u$ and its directed neighbours $\{y_i \in \vec{\mathbf{y}} : u \sim y_i\}$ in $\vec{\mathbf{y}}$.

In particular, $\lambda_{n(\vec{y})+1}(\cdot \mid \vec{y})$ is invariant under permutations of $\vec{y}$. In the sequel, we shall sometimes abuse notation somewhat and write $\lambda_{n(\mathbf{y})+1}(\cdot \mid \mathbf{y})$, where $\mathbf{y} = \{y_i \in \vec{y}\}$, and $n(\mathbf{y})$ is the cardinality of the set $\mathbf{y}$.

Due to the lack of symmetry, sequential spatial processes are particularly useful for modelling inhomogeneous space-time processes. For instance, to obtain a complete packing of non-intersecting particles [28], a typical algorithm keeps adding particles randomly until this becomes impossible due to violation of the condition of no overlap, see e.g. the recent review papers [6, 27]. A related soft-core model for landslides was proposed by Fiksel and Stoyan [9].

**Example 1 (Simple sequential inhibition).** Let $\pi(\cdot)$ be a probability density with respect to Lebesgue measure on a bounded planar Borel set $D$ of strictly positive area. Suppose a population of animals arrives in the region $D$ to set up nests, and would like to pick locations according to $\pi(\cdot)$. In other words, those subregions that have high $\pi$-mass are deemed the most suitable. However, as each animal needs its space, a newly arriving animal may not build its nest closer than some $r > 0$ to an existing nest location [4, 5]. Thus, if we write $d(z, \vec{x}) = \min\{||z - x|| : x \in \vec{x}\}$ for the minimal distance between $z \in D$ and the components of $\vec{x}$,

$$p_n(x_1, \ldots, x_n) \propto \pi(x_1) \frac{\pi(x_2) \mathbf{1}\{d(x_2, x_1) > r\}}{\int_D \pi(z) \mathbf{1}\{d(z, x_1) > r\} \, dz} \cdots \frac{\pi(x_n) \mathbf{1}\{d(x_n, \vec{x}_{<n}) > r\}}{\int_D \pi(z) \mathbf{1}\{d(z, \vec{x}_{<n}) > r\} \, dz} \quad (9)$$

for $x_i \in D$, $i = 1, \ldots, n$. Some care has to be taken about division by zero. If $n$ is larger than the packing number for balls of radius $r$ in $D$, or $q_n = 0$ for some other reason, $p_n(\cdot, \ldots, \cdot)$ may be chosen to be any probability density on $D^n$. For smaller $n$, or if $q_n > 0$, set $p_n(x_1, \ldots, x_n) = 0$ whenever some term in the denominator of (9) is zero, and renormalise $p_n(\cdot, \ldots, \cdot)$ to integrate to unity. The total number of animals in $D$ is either fixed – as in the formulation by [4] – or random according to some probability mass function $q_n$, $n \leq n_p$, the packing number, as in *Random Sequential Adsorption* [6, 27].

If $n$ is fixed, that is $q_n = \mathbf{1}\{n = n_0\}$ for some $n_0 \in \mathbb{N}$, $f(\vec{x})$ is not hereditary, therefore not Markovian in the sense of Definition 1. If $q_n > 0$ for $n \leq n_p$, and $\sum_{n=0}^{n_p} q_n = 1$, the density is hereditary. In that case, provided $f(x_1, \ldots, x_n) > 0$,

$$\lambda_{n+1}(u \mid (x_1, \ldots, x_n)) = \frac{c_{n+1} \, q_{n+1}}{c_n \, q_n} \frac{\pi(u) \mathbf{1}\{d(u, \vec{x}) > r\}}{\int_D \pi(z) \mathbf{1}\{d(z, \vec{x}) > r\} \, dz}$$

(with $0/0 = 0$), where $c_n$ is the normalising constant of $p_n(\cdot, \ldots, \cdot)$. The likelihood ratio is invariant under permutations of the components of $\vec{x}$. However, it may depend on the sequence length through $c_{n+1} q_{n+1}/(c_n q_n)$, and, moreover, $\int_D \pi(z) \mathbf{1}\{d(z, \vec{x}) > r\} \, dz$ may depend on the geometry of $\vec{x}$, i.e. on the whole sequence. Hence, in general, $f(\cdot)$ is Markovian only with respect to the trivial relation in which each pair of points is related.

**Example 2 (Sequential soft core).** Further to Example 1, suppose the animals have no locational preferences, but claim territory within a certain radius. The radius can depend deterministically on the location [14] (in fertile regions, less space is needed than in poorer ones), be random and captured by a mark, or a combination of both [1]. To be specific, suppose an animal settling at $x \in D$ claims

---

[1] Assume the mark space $(0, \infty)$, for concreteness, is equipped with a density $g(\cdot)$. By the transformation $h_{(c_1, c_2)}(x, m) = (c_1 x, c_2 m)$ for $c_1, c_2 > 0$, the Lebesgue intensity measure $\mu(\cdot)$ is scaled

the region within a (stochastic) radius $m > 0$ of $x$ for itself, and that newcomers are persuaded to avoid this region. An appropriate density with respect to $\nu(\cdot)$ could be

$$f(\vec{\mathbf{y}}) \propto \beta^{n(\vec{\mathbf{y}})} \prod_{i=1}^{n(\vec{\mathbf{y}})} \gamma^{\sum_{j<i} \mathbf{1}\{||x_i - x_j|| \leq m_j\}},$$

$\vec{\mathbf{y}} \ni y_i = (x_i, m_i) \in D \times (0, \infty)$. Here $\beta > 0$ is an intensity parameter, and $0 \leq \gamma < 1$ reflects the strength of persuasion, the smaller $\gamma$, the stronger the inhibition. Indeed, for $\gamma = 0$, no new arrival is allowed to enter the territory of well-settled animals. Alternatively, if invaders demand space according to their own mark, replace $m_j$ by $m_i$ in the exponent of $\gamma$.

Note that

$$(n(\vec{\mathbf{y}}) + 1)\, \lambda_{n(\vec{\mathbf{y}})+1}((x, m) \mid \vec{\mathbf{y}}) = \beta\, \gamma^{\sum_{j=1}^{n(\vec{\mathbf{y}})} \mathbf{1}\{||x - x_j|| \leq m_j\}}$$

depends only on those $y_j = (x_j, m_j)$ for which $||x - x_j|| \leq m_j$. Hence $f(\cdot)$ is Markov with respect to the relation $(x, r) \sim (y, s) \Leftrightarrow ||x - y|| \leq s$.

## 3. Hammersley–Clifford factorisation

The goal of this section is to formulate and prove a factorisation theorem for sequential spatial processes in analogy with (2) and (6).

To do so, we need the notion of a directed clique interaction function. Indeed, for a marked point $z \in D \times M$, define the sequence $\vec{\mathbf{y}}$ to be a $z$-clique with respect to a reflexive relation $\sim$ on $D \times M$ if $\vec{\mathbf{y}}$ either has length zero or all its components $y \in \vec{\mathbf{y}}$ satisfy $z \sim y$. The definition is $z$-directed but otherwise permutation invariant, so we may map $\vec{\mathbf{y}}$ onto the set $\mathbf{y} \in \widetilde{N^{\mathrm{f}}}$, the family of unordered finite marked point configurations, by ignoring the permutation.

Our main theorem is the following.

**Theorem 1.** *A sequential spatial process with density $f(\cdot)$ is Markov with respect to $\sim$ if and only if it can be factorised as*

$$f(y_1, \ldots, y_n) = f(\emptyset) \prod_{i=1}^{n} \prod_{\mathbf{z} \subseteq \mathbf{y}_{<i}} \varphi(y_i, \mathbf{z}) \tag{10}$$

*for some non-negative, jointly measurable interaction function $\varphi : (D \times M) \times \widetilde{N^{\mathrm{f}}} \to [0, \infty)$ that vanishes except on cliques (i.e. $\varphi(u, \mathbf{z}) = 1$ if $\mathbf{z}$ is no $u$-clique with respect to $\sim$). Here $\mathbf{y}_{<i} = \{y_1, \ldots, y_{i-1}\}$.*

---

by $c_1^{-2}$ on $c_1 D$, the new mark density is $c_2^{-1} g(m/c_2)$. Let $c(\cdot, \cdot)$ be a scaling function with component functions $c_i(\cdot, \cdot)$, $i = 1, 2$, such that $0 < \underline{c} \leq c_i(x, m) \leq \overline{c} < \infty$ uniformly. Let $Y$ be a template marked point process with hereditary density $f(\cdot)$, conditional intensity $\lambda(\cdot \mid \cdot)$, and define a sequential conditional intensity with respect to $d\mu_c(x, m) = c_1(x, m)^{-2} c_2(x, m)^{-1} g(m/c_2(x, m))$ by $\lambda_{Y_c}((x, m) \mid \mathbf{y}) := \lambda\left(\left(\frac{x}{c_1(x,m)}, \frac{m}{c_2(x,m)}\right) \mid h^{-1}_{c(x,m)}(\mathbf{y})\right) / (n(\mathbf{y}) + 1)$. Provided $f(\cdot)$ is locally stable, the sequential density $f_{Y_c}(\cdot)$ defined by $\lambda_{Y_c}(\cdot \mid \cdot)$ as in (8) is well-defined. For example, the marked hard core density $f(\mathbf{y}) \propto \prod \mathbf{1}\{||x_i - x_j|| > m_i + m_j\}$ with $\lambda((x, m) \mid \mathbf{y}) = \mathbf{1}\{||x - x_i|| > m + m_i \text{ for all } (x_i, m_i) \in \mathbf{y}\}$ is transformed into $f_{Y_c}(\vec{\mathbf{y}}) \propto \prod_{i<j} \mathbf{1}\left\{||x_i - x_j|| > \frac{c_1(x_j, m_j)}{c_2(x_j, m_j)}(m_i + m_j)\right\}$. Intuitively, the resulting process $Y_c$ looks like a scaling of $Y$ at $(x, m)$ by $c(x, m)$. Indeed, if $Y$ is Markov with respect to $\sim$, then $Y_c$ is Markov with respect to $(u, r) \sim_c (v, s) \Leftrightarrow \left(\frac{u}{c_1(u,r)}, \frac{r}{c_2(u,r)}\right) \sim \left(\frac{v}{c_1(u,r)}, \frac{s}{c_2(u,r)}\right)$ and the $\sim_c$-neighbourhood $\partial_c(x, m) = h_{c(x,m)}(\partial(h^{-1}_{c(x,m)}(x, m)))$ inherits from $\partial(\cdot)$ geometric properties such as convexity that are invariant under rescaling.

*Proof.* To show that any density of the form (10) is Markovian, suppose $f(\vec{y}) > 0$ but $f(\vec{z}) = 0$ for some subsequence $\vec{z}$ of $\vec{y}$. Then, again with the notation $n(\vec{z})$ for the length of $\vec{z}$,

$$f(\emptyset) \prod_{i=1}^{n(\vec{z})} \prod_{\mathbf{x} \subseteq \mathbf{z}_{<i}} \varphi(z_i, \mathbf{x}) = 0.$$

Hence either $f(\emptyset) = 0$ or some term $\varphi(z_i, \mathbf{x}) = 0$, but, as both feature in the factorisation (10) of $f(\vec{y})$ too, $f(\vec{y})$ must be 0 as well, in contradiction with the assumption. Assume that $f(\vec{y}) > 0$. Now, $f((\vec{y}, u))/f(\vec{y}) = \prod_{\mathbf{z} \subseteq \mathbf{y}} \varphi(u, \mathbf{z})$ depends only on $u$ and its directed neighbours in $\mathbf{y}$, as $\varphi(u, \mathbf{z}) = 1$ whenever $\mathbf{z}$ is no $u$-clique.

Reversely, set $\varphi(y_1, \emptyset) := f(y_1)/f(\emptyset) = \lambda_1(y_1|\emptyset)$, and define putative interaction functions $\varphi(\cdot, \cdot)$ recursively as follows. If $\{y_1, \ldots, y_n\}$ is no $y_{n+1}$-clique, $\varphi(y_{n+1}, \{y_1, \ldots, y_n\}) := 1$; else

$$\varphi(y_{n+1}, \{y_1, \ldots, y_n\}) := \frac{f(y_1, \ldots, y_{n+1})}{f(y_1, \ldots, y_n) \prod_{\mathbf{z}} \varphi(y_{n+1}, \mathbf{z})} \tag{11}$$

where the product ranges over all strict subsets $\{y_1, \ldots, y_n\} \neq \mathbf{z} \subset \{y_1, \ldots, y_n\}$. To deal with any zeroes, we use the convention $0/0 = 0$.

We first show that $\varphi(\cdot, \cdot)$ is a well-defined interaction function. By the Markov assumption, $f(y_1, \ldots, y_{n+1})/f(y_1, \ldots, y_n)$ is invariant under permutations of $y_1, \ldots, y_n$, so a simple induction argument yields the permutation invariance of $\varphi(y_{n+1}, \{y_1, \ldots, y_n\})$ in its second argument. By definition, $\varphi(\cdot, \cdot)$ vanishes except on cliques, hence it is an interaction function. To show that if the denominator of (11) is zero, so is the numerator, note that if $f(\emptyset) = 0$, by the assumption that the process is hereditary, necessarily $f \equiv 0$ which contradicts the fact that $f(\cdot)$ is a probability density. Therefore $f(\emptyset) > 0$ and $\varphi(y_1, \emptyset)$ is well-defined. Suppose $\varphi(\cdot, \cdot)$ is well-defined for sets of cardinality at most $n - 1 \geq 0$ as its second argument. Let $\{y_1, \ldots, y_n\}$ be an $y_{n+1}$-clique. The Markov assumption implies that if $f(y_1, \ldots, y_n) = 0$, also $f(y_1, \ldots, y_{n+1})$ is zero. Furthermore, if $f(y_1, \ldots, y_n) > 0$ but $\varphi(y_{n+1}, \mathbf{z}) = 0$ for some strict subset $\mathbf{z} \subset \{y_1, \ldots, y_n\}$, by the induction assumption, $f((\vec{z}, y_{n+1})) = 0$, where the sequence $\vec{z}$ is obtained from $\mathbf{z}$ by the permutation induced by $(y_1, \ldots, y_n)$. A fortiori $f(y_1, \ldots, y_{n+1}) = 0$. We conclude that (11) is well-defined by induction.

It remains to show that $\varphi(\cdot, \cdot)$ satisfies the desired factorisation. To do so, we again proceed by induction. By definition, $f(y_1) = f(\emptyset) \varphi(y_1, \emptyset)$ for any $y_1 \in D \times M$, so the factorisation holds for sequences of length at most 1. Suppose (10) holds for all sequences that are at most $n \geq 1$ long, and consider any sequence $\vec{y} = (y_1, \ldots, y_{n+1})$ with components in $D \times M$. If $f(y_1, \ldots, y_n) = 0$, by the assumption on hereditariness, $f(y_1, \ldots, y_{n+1}) = 0$. By the induction hypothesis, $f(y_1, \ldots, y_n) = f(\emptyset) \prod_{i=1}^{n} \prod_{\mathbf{z} \subseteq \mathbf{y}_{<i}} \varphi(y_i, \mathbf{z}) = 0$ which implies $f(\emptyset) \prod_{i=1}^{n+1} \prod_{\mathbf{z} \subseteq \mathbf{y}_{<i}} \varphi(y_i, \mathbf{z}) = 0 = f(y_1, \ldots, y_{n+1})$. Hence without loss of generality we may assume that $f(y_1, \ldots, y_n)$ is strictly positive. Then

$$\begin{aligned} f(y_1, \ldots, y_{n+1}) &= \frac{f(y_1, \ldots, y_{n+1})}{f(y_1, \ldots, y_n)} f(y_1, \ldots, y_n) \\ &= \frac{f(y_1, \ldots, y_{n+1})}{f(y_1, \ldots, y_n)} f(\emptyset) \prod_{i=1}^{n} \prod_{\mathbf{z} \subseteq \mathbf{y}_{<i}} \varphi(y_i, \mathbf{z}) \end{aligned} \tag{12}$$

because of the induction hypothesis. We shall distinguish several cases.

Firstly, suppose that $\{y_1, \ldots, y_n\}$ is an $y_{n+1}$-clique and recall that if the denominator in the right hand side of (11) is zero, then so are $\varphi(y_{n+1}, \{y_1, \ldots, y_n\})$ and $f(y_1, \ldots, y_{n+1})$, so the desired factorisation holds. If, on the other hand, the denominator is strictly positive, by definition (11),

$$\frac{f(y_1, \ldots, y_{n+1})}{f(y_1, \ldots, y_n)} = \prod_{\mathbf{z} \subseteq \{y_1, \ldots, y_n\}} \varphi(y_{n+1}, \mathbf{z}) \tag{13}$$

and (10) holds by substitution of (13) in (12).

Secondly, assume $\{y_1, \ldots, y_n\}$ is no $y_{n+1}$-clique. Then, $y_{n+1} \not\sim y_i$ for some $y_i$, $i = 1, \ldots, n$. Now, by the Markov assumption, $f(y_1, \ldots, y_{i-1}, y_{i+1}, \ldots, y_n) > 0$, and

$$\begin{aligned}
\lambda_{n+1}(y_{n+1} \mid \{y_1, \ldots, y_n\}) &= \lambda_{n+1}(y_{n+1} \mid \{y_1, \ldots, y_{i-1}, y_{i+1}, \ldots, y_n\}) \\
&= \frac{1}{n+1} \frac{f(y_1, \ldots, y_{i-1}, y_{i+1}, \ldots, y_{n+1})}{f(y_1, \ldots, y_{i-1}, y_{i+1}, \ldots, y_n)} \\
&= \frac{1}{n+1} \prod_{\mathbf{z} \subseteq \{y_1, \ldots, y_{i-1}, y_{i+1}, \ldots, y_n\}} \varphi(y_{n+1}, \mathbf{z}) \\
&= \frac{1}{n+1} \prod_{\mathbf{z} \subseteq \{y_1, \ldots, y_n\}} \varphi(y_{n+1}, \mathbf{z})
\end{aligned}$$

where the last two equations follow from the induction hypothesis and the fact that $\varphi(y_{n+1}, \mathbf{z}) = 1$ whenever $\mathbf{z}$ contains $y_i$. We conclude that (13) holds, hence (10). □

The factorisation (10) is similar to that of the joint distribution of a Markov chain, cf. (2). Indeed, the product of transition probabilities $p(x_{i-1}, x_i)$ from state $x_{i-1}$ to $x_i$ is replaced by a product of interaction functions $\prod_{\mathbf{z} \subseteq \mathbf{y}_{<i}} \varphi(y_i, \mathbf{z})$ over cliques of 'past' neighbours of $y_i$ similar to (6).

Theorem 1 provides a third way of defining a sequential process: by specifying its clique interaction functions. Of course one must verify integrability for each choice.

**Example 1 (ctd).** Recall that we need to assume $q_n > 0$ for $n \leq n_p$ and $q_n = 0$ for $n > n_p$ in order to be sure that the simple sequential inhibition model (9) is hereditary.

Under this assumption, the interaction functions can be computed iteratively by (11). Set $r_n := n c_n q_n / (c_{n-1} q_{n-1})$ for $n = 1, 2, \ldots, n_p$, 0 otherwise, and $I(\mathbf{x}) := \int_D \pi(z) \mathbf{1} \{d(z, \mathbf{x}) > r\} \, dz$ for any finite set $\mathbf{x}$ of points in $D$, with $I(\emptyset) = 1$. Then

$$\varphi(x, \emptyset) = r_1 \pi(x);$$

$$\varphi(x, \mathbf{y}) = \mathbf{1} \{d(x, \mathbf{y}) > r\} \exp \left[ \sum_{\mathbf{z} \subseteq \mathbf{y}} (-1)^{n(\mathbf{y} \setminus \mathbf{z})} \log \left( \frac{r_{n(\mathbf{z})+1}}{I(\mathbf{z})} \right) \right]$$

for non-empty configurations $\mathbf{y}$ with $I(\mathbf{y}) > 0$. Otherwise, $\varphi(x, \mathbf{y}) = 0$.

In order to verify the above expressions, note that by equation (8), it is sufficient to verify that the sequential conditional intensity, or equivalently the likelihood ratio for adding a point at the end of a given sequence, has the desired form. Indeed, for $\vec{\mathbf{x}} = (x_1, \ldots, x_n)$, $n \geq 1$, such that $I(\mathbf{x})$ and $f(\vec{\mathbf{x}})$ are strictly positive,

$$\begin{aligned}
\frac{f(s_{n+1}(\vec{\mathbf{x}}, x_{n+1}))}{f(\vec{\mathbf{x}})} &= r_1 \pi(x_{n+1}) \exp \left[ \sum_{\emptyset \neq \mathbf{y} \subseteq \mathbf{x}} \sum_{\mathbf{z} \subseteq \mathbf{y}} (-1)^{n(\mathbf{y} \setminus \mathbf{z})} \log \left( \frac{r_{n(\mathbf{z})+1}}{I(\mathbf{z})} \right) \right] \\
&= \frac{r_{n+1} \pi(x_{n+1})}{I(\mathbf{x})}
\end{aligned}$$

provided $d(x_{n+1}, \{x_1, \ldots, x_n\}) > r$, which implies $d(x_{n+1}, \mathbf{y}) > r$ for any subset $\mathbf{y}$ of $\mathbf{x} = \{x_1, \ldots, x_n\}$, and zero otherwise. The last identity is a consequence of Newton's binomium.

**Example 2 (ctd).** For the sequential soft core model,

$$\begin{aligned}
\varphi(y, \emptyset) &= \beta; \\
\varphi((x_1, m_1), (x_2, m_2)) &= \gamma^{\mathbf{1}\{||x_1 - x_2|| \leq m_2\}},
\end{aligned}$$

and 1 otherwise [2].

**Example 3.** Pairwise interaction models of the form

$$f(\emptyset) \prod_{i=1}^{n} \left[ \varphi(y_i, \emptyset) \prod_{j<i} \varphi(y_i, y_j) \right]$$

are particularly convenient. Note that the sequential soft core model discussed in the previous example exhibits pairwise interactions only. It has a constant penalty term $\gamma \in [0, 1)$ for each pair $y_i \sim y_j$. It may be more natural to let the interaction vary with distance, for example quadratically as in

$$\varphi((x, r), (y, s)) = 1 - (1 - ||x - y||^2 / R_{r,s}^2)^2$$

for $||x - y|| \leq R_{r,s}$. Here the range of interaction $R_{r,s}$ may depend on the marks. A sufficient condition for the model to be integrable is that $\varphi(y, \emptyset)$ is uniformly bounded in $y$ and the pairwise interaction function $\varphi(y_i, y_j)$ is bounded by 1. Thus, pairwise interaction models are particularly suitable for modelling repulsion between marked points.

## 4. Dynamic representation

Markov sequential spatial processes arise naturally as the limit distribution of a jump process [7] with transitions that add or remove one component of the sequence at a time. From a statistical point of view, it is then possible to obtain samples from some sequential spatial model of interest by running the pure jump process into equilibrium. Standard Markov chain Monte Carlo ideas [10] can then be applied to estimate model parameters by the maximum likelihood principle, to perform likelihood ratio tests to assess the goodness-of-fit of the model, to compute confidence regions, profile likelihoods and so on. Moreover, at least in principle, it is possible to determine exactly when equilibrium is reached, by applying the coupling ideas of Propp and Wilson [24] along the lines proposed in [8, 16, 19] and the references therein.

Below, we first propose a birth-and-death jump process, then develop a discrete time Metropolis–Hastings style sampler.

---

[2] The density of the scaled process $Y_c$ of Example 2 can be factorised as $f_{Y_c}(\emptyset) \prod_{i=1}^{n(\mathbf{y})} \prod_{\mathbf{z} \subseteq \mathbf{y}_{<i}} \varphi \left( h_{c(y_i)}^{-1}(\mathbf{z}) \cup \{(\frac{x_i}{c_1(y_i)}, \frac{m_i}{c_2(y_i)})\} \right)$ if $\varphi(\cdot)$ denotes the interaction function of a template Markov marked point process $Y$.

### 4.1. Spatial birth-and-death processes

A *spatial birth-and-death process* is a continuous time Markov process with state space $N^{\mathrm{f}}$. Its only transitions are the insertion of a marked point in the current sequence (a *birth*), or the deletion of a component (a *death*). Suppose the current state is $\vec{\mathbf{y}}$, and write $b_i(\vec{\mathbf{y}}, u) \, d\mu \times \mu_M(u)$ for the birth rate of a new marked point in $du$, $u \in D \times M$, to be inserted at position $i$ of $\vec{\mathbf{y}}$, and $d_i(\vec{\mathbf{y}})$ for the death rate of the reverse transition. Then the detailed balance equations are given by

$$\frac{e^{-\mu(D)}}{n!} f(\vec{\mathbf{y}}) \, b_i(\vec{\mathbf{y}}, u) = \frac{e^{-\mu(D)}}{(n+1)!} f(s_i(\vec{\mathbf{y}}, u)) \, d_i(s_i(\vec{\mathbf{y}}, u)).$$

If we take unit death rate for each marked point,

$$b_i(\vec{\mathbf{y}}, u) = \frac{f(s_i(\vec{\mathbf{y}}, u))}{(n+1) \, f(\vec{\mathbf{y}})} \qquad (14)$$

again with the convention $0/0 = 0$. Note that if the jump process starts in the set $\{\vec{\mathbf{y}} \in N^{\mathrm{f}} : f(\vec{\mathbf{y}}) > 0\}$, it will almost surely never leave it. Clearly, for this choice of birth and death rates, $f(\cdot)$ is an invariant probability density.

Since our process evolves in continuous time, some care has to be taken to avoid explosion, that is, infinitely many transitions in finite time. Set

$$\mathcal{B}_n := \sup_{n(\vec{\mathbf{y}}) = n} \sum_{i=1}^{n+1} \int_{D \times M} b_i(\vec{\mathbf{y}}, u) \, d\mu \times \mu_M(u)$$

for the upper bound on the total birth rate from sequences of length $n \geq 0$, and note that

$$\delta_n := \inf_{n(\vec{\mathbf{y}}) = n} \sum_{i=1}^{n} d_i(\vec{\mathbf{y}}) = n.$$

By [23, Prop. 5.1, Thm. 7.1], sufficient conditions for the existence of a unique jump process with birth rates (14) and unit death rates with a unique invariant probability measure to which it converges in distribution regardless of the initial state are that $\delta_n > 0$ for all $n \geq 1$ and one of the following holds:

- $\mathcal{B}_n = 0$ for all sufficiently large $n \geq n_0 \geq 0$;
- $\mathcal{B}_n > 0$ for all $n \geq 1$ and

$$\sum_{n=1}^{\infty} \frac{\delta_1 \cdots \delta_n}{\mathcal{B}_1 \cdots \mathcal{B}_n} = \infty; \qquad \sum_{n=2}^{\infty} \frac{\mathcal{B}_1 \cdots \mathcal{B}_{n-1}}{\delta_1 \cdots \delta_n} < \infty.$$

For densities $f(\cdot)$, such as the sequential soft core model of Example 2, that satisfy

$$f(s_i(\vec{\mathbf{y}}, u)) \leq \beta f(\vec{\mathbf{y}}) \qquad (15)$$

for some $\beta > 0$ and all $\vec{\mathbf{y}} \in N^{\mathrm{f}}$, $u \in D \times M$, $i = 1, \ldots, n(\vec{\mathbf{y}}) + 1$, the total birth rate is bounded by $\beta \mu(D)$, and it is easily verified that the above Preston conditions hold. Moreover, from a practical point of view, the jump process may be implemented by thinning the transitions of a process with unit death rate and constant birth rate $\tilde{b}_i(\vec{\mathbf{y}}, \cdot) \equiv \beta/(n(\vec{\mathbf{y}}) + 1)$, which avoids having to compute explicitly the parameter

$$\sum_{i=1}^{n(\vec{\mathbf{y}})+1} \int_{D \times M} b_i(\vec{\mathbf{y}}, u) \, d\mu \times \mu_M(u) + n(\vec{\mathbf{y}})$$

of the exponentially distributed waiting times in between jumps. The retention probability of a transition from $\vec{y}$ to $s_i(\vec{y}, u)$ is given by

$$\frac{f(s_i(\vec{y}, u))}{\beta f(\vec{y})} = \frac{1}{\beta} \prod_{\mathbf{z} \subseteq \mathbf{y}_{<i}} \varphi(u, \mathbf{z}) \prod_{j=i}^{n(\vec{y})} \prod_{\mathbf{z} \subseteq \mathbf{y}_{<j}} \varphi(y_j, \mathbf{z} \cup \{u\}).$$

The first term depends on the directed neighbours $y_j$ of $u$ ($u \sim y_j$) in $\mathbf{y}_{<i}$. The interaction functions in the last product reduce to 1 if $y_j \not\sim u$.

### 4.2. A Metropolis–Hastings sampler

A popular, flexible and generally applicable technique for generating samples from a complex or high-dimensional probability density is the *Metropolis–Hastings method*, see e.g. [10]. Briefly, the method works by proposing an update according to a distribution that is convenient to sample from, and then to accept or reject this proposal with a probability that is chosen so as to make sure the detailed balance equations are satisfied.

In the context of this paper, two types of proposals are considered: births and deaths. More precisely, given the current state is $\vec{y}$, with probability $1/2$ propose a birth, otherwise a death. In the first case, select a position $i$ to insert and a new marked point $u$ uniformly on $1, \ldots, n(\vec{y}) + 1$ respectively $D \times M$ (w.r.t. $\mu \times \mu_M$), and accept the new state $s_i(\vec{y}, u)$ with probability

$$\alpha(\vec{y}, s_i(\vec{y}, u)) := \min\left\{1, \frac{f(s_i(\vec{y}, u)) \, \mu(D)}{f(\vec{y}) \, (n(\vec{y}) + 1)}\right\}.$$

In case a death is proposed, if $\vec{y}$ is the empty sequence, do nothing; otherwise, select the position of the marked point to be removed uniformly, say $i$, and accept the transition with probability

$$\alpha(\vec{y}, \vec{y}_{(-i)}) := \min\left\{1, \frac{f(\vec{y}_{(-i)}) \, n(\vec{y})}{f(\vec{y}) \, \mu(D)}\right\}$$

where $\vec{y}_{(-i)}$ denotes the subsequence of $\vec{y}$ obtained by removing the $i^{\text{th}}$ component. It is easily seen that $f(\cdot)$ is an invariant density. If we start the chain in a sequence $\vec{y}$ for which $f(\vec{y}) > 0$, the chain will almost surely never leave the set $H := \{\vec{y} \in N^{\text{f}} : f(\vec{y}) > 0\}$. We shall restrict the Markov chain to $H$ from now on.

In order to show that the Metropolis–Hastings chain $Y_n$, $n \in \mathbb{N}_0$, converges to $f(\cdot)$ in total variation from any initial state in $H$, we need to establish aperiodicity and Harris recurrence, that is *Harris ergodicity* [21]. Aperiodicity immediately follows from the fact that self-transitions occur. A sufficient condition for Harris recurrence (in fact for the even stronger property of geometric ergodicity) is that there exist a function $V : N^{\text{f}} \cap H \to [1, \infty)$, constants $b < \infty$ and $\gamma < 1$, and a measurable small set $C \subset N^{\text{f}} \cap H$ such that

$$PV(\vec{y}) := \mathbb{E}\left[V(X_{n+1})|X_n = \vec{y}\right] \leq \gamma V(\vec{y}) + b \, \mathbf{1}\{\vec{y} \in C\}. \tag{16}$$

Recall that a set $C$ is small if there exists a non-zero measure $\phi(\cdot)$ and an integer $n$ such that the probability that the Metropolis–Hastings chain reaches any measurable subset $B$ of $N^{\text{f}}$ from any $\vec{y} \in C$ in $n$ steps is at least as large as $\nu(B)$. See [21] for further details.

From now on, assume that (15) holds. Then, the drift condition (16) can be verified by an adaptation of the proof of [12, Prop. 3.3] with $V(\vec{y}) = A^{n(\vec{y})}$ for some $A > 1$. Indeed, note that the acceptance probability for inserting $u \in D \times M$ at position $i$ in $\vec{y}$ is bounded by $\beta \mu(D)/(n(\vec{y})+1)$, which does not exceed a prefixed constant $\epsilon > 0$ if $n(\vec{y})$ is sufficiently large. Similarly, the acceptance probability for removing a component from the sequence $\vec{y}$ reduces to 1 if $\vec{y}$ is long enough. Now,

$$PV(\vec{y}) = \frac{1}{2} V(\vec{y}) \sum_{i=1}^{n(\vec{y})+1} \int_{D \times M} \frac{A-1}{\mu(D)(n(\vec{y})+1)} \alpha(\vec{y}, s_i(\vec{y}, u)) \, d\mu \times \mu_M(u)$$
$$+ \frac{1}{2} V(\vec{y}) \sum_{i=1}^{n(\vec{y})} \frac{A^{-1}-1}{n(\vec{y})} \alpha(\vec{y}, \vec{y}_{(-i)}) + V(\vec{y}). \tag{17}$$

Since for sufficiently long sequences of marked points, say $n(\vec{y}) > N_\epsilon$, deaths are always accepted and $\alpha(\vec{y}, s_i(\vec{y}, u)) \leq \epsilon$ uniformly in its arguments, (17) is less than or equal to $\left[\frac{1}{2}(A-1)\epsilon + \frac{1}{2}(A^{-1}-1) + 1\right] V(\vec{y})$. This and the uniform lower bound $1/(\beta \mu(D))$ on the acceptance probability for deaths yield the desired result by the same arguments as in the proof of Propositions 3.2–3.3 in [12] with $C = \cup_{n=0}^{N_\epsilon} (D \times M)^n \cap H$.

### 4.3. Reversible jump Markov chains

The Metropolis–Hastings algorithm of the previous subsection is a special case of a so-called reversible jump Markov chain [13]. This is a proposal–acceptance/rejection scheme for moves between subspaces of different dimension. Suppose a target equilibrium distribution and a proposal probability for each move type are given. The aim is to define acceptance probabilities in such a way that the resulting Markov chain is well-defined and the detailed balance equations hold. In order to avoid singularities due to the different dimensions, for each move type, a symmetric dominating measure $\xi(\cdot, \cdot)$ on the product of the state space with itself is needed. It turns out [13] that acceptance probabilities that are ratios of the joint density of the target and proposal distributions with respect to $\xi(\cdot, \cdot)$ evaluated at the proposed and current state do the trick.

In the context of Section 4.2, define, for measurable subsets $A$ and $B$ of $N^{\mathrm{f}}$,

$$\xi(A \times B) = \int_{N^{\mathrm{f}}} \mathbf{1}\{\vec{y} \in A\} \frac{1}{n(\vec{y})+1}$$
$$\times \sum_{i=1}^{n(\vec{y})+1} \int_{D \times M} \mathbf{1}\{s_i(\vec{y}, u) \in B\} \, d\mu \times \mu_M(u) \, d\nu(\vec{y})$$
$$+ \int_{N^{\mathrm{f}}} \mathbf{1}\{\vec{y} \in A\} \sum_{i=1}^{n(\vec{y})} \mathbf{1}\{\vec{y}_{(-i)} \in B\} \, d\nu(\vec{y}).$$

The measure is symmetric by the form of the reference measure $\nu(\cdot)$.

The joint bi-variate density of the product of the distribution to sample from and the birth proposal distribution with respect to $\xi(\cdot, \cdot)$ is given by $f(\vec{y}, s_i(\vec{y}, u)) = f(\vec{y})/(2\mu(D))$. Similarly for death proposals we have $f(\vec{y}, \vec{y}_{(-i)}) = f(\vec{y})/(2n(\vec{y}))$. Hence, the acceptance probabilities for births and deaths, as given in the previous section, follow as ratios of joint bi-variate densities truncated at 1. Note that there is no need to worry about division by zero, as we had restricted the Metropolis–Hastings chain to the family of sequences having strictly positive density.

**Acknowledgments**

The author would like to thank Mike Keane for piloting CWI's spatial statistics research safely through a re-organisation, for luring her back to The Netherlands, and for a pleasant collaboration on the annual Lunteren Meeting. She acknowledges A. Baddeley, U. Hahn, E. Vedel Jensen, J. Møller, L. Nielsen, Y. Ogata, T. Schreiber, R. Stoica, D. Stoyan, and the participants in STW project CWI.6156 for discussions and preprints.